\documentclass[12pt,letterpaper]{article}
\usepackage{amsfonts}
\usepackage{amssymb}
\usepackage{latexsym}  
\usepackage{pifont}  
\usepackage{amsmath} 
\usepackage{bm} 

\newcommand\DD{{\mathcal D}}
\newcommand\FF{{\mathcal F}}

\newcommand\duX{{\widehat X}}
\newcommand\wCC{\widetilde{\CC}}
\newcommand\wDD{\widetilde{\DD}}
\newcommand\CC{{\mathcal C}}

\newcommand\vphi{\boldsymbol{\varphi}}  

\newcommand\Fphi{\FF_{\vphi}}

\newcommand\CCC{{\mathbb C}}
\newcommand\TTT{{\mathbb T}}
\newcommand\ZZZ{{\mathbb Z}}

\newcommand\tttt{\mathfrak{t}}
\newcommand\gggg{\mathfrak{g}}

\newcommand\onto{\twoheadrightarrow}  

\newcommand\res{\mathord {\upharpoonright}}  

\newcommand\Hom{\mathrm{Hom}}  
\newcommand\cl{\mathrm{cl}}    

\newcommand\iv{^{-1}} 

\newcommand\one{\raisebox{-0.2ex}{\mbox{\large \ding{"AC}}}}
\newcommand\bohr{\mathsf{b}}

\newcommand\cat{^{\mathord{\frown}}}  

\newcommand\eop{$\ \ {\vcenter
   {\hrule
   \hbox{\vrule height 9pt \kern 9pt \vrule height 9pt}
   \hrule}}$\vskip 1.0 pt}

\newenvironment{itemizz}{\begin{itemize}\setlength{\itemsep}{-1mm}} %
{\end{itemize}}                              

\newenvironment{itemizn}[1] 
{\begin{itemize} \setlength{\itemsep}{-1mm} %
\renewcommand\labelitemi{\ding{#1}}} %
{\end{itemize}}

\newtheorem{theorem}{Theorem}[section]
\newtheorem{definition}[theorem]{Definition}
\newtheorem{lemma}[theorem]{Lemma}

\newtheorem{proposition}[theorem]{Proposition}
\newtheorem{example}[theorem]{Example}

\newenvironment{proof}{{\bf Proof.}}{\eop\medskip}
\newenvironment{proofof}[1]{\medskip \textbf{Proof of #1.}}{\eop\medskip}

\begin{document}

\title{
Limits in Compact Abelian Groups\footnote{
2000 Mathematics Subject Classification:
Primary 54H11, 22C05; Secondary  43A40.
Key Words and Phrases: Compact group, character, Bohr topology,
pointwise convergence.
}}

\author{Joan E. Hart\footnote{University of Wisconsin, Oshkosh,
WI 54901, U.S.A.,
\ \ hartj@uwosh.edu}
\  and
Kenneth Kunen\footnote{University of Wisconsin,  Madison, WI  53706, U.S.A.,
\ \ kunen@math.wisc.edu}
\thanks{Both authors partially supported by NSF Grant DMS-0097881.}
}

\maketitle

\begin{abstract}
For $X$ a compact abelian group and
$B$ an infinite subset of its dual $\widehat X$, 
let $\CC_B$ be the set of all $x \in X$ such that 
$\langle \varphi(x) : \varphi \in B\rangle$ converges to $1$.
If $\FF$ is a free filter on $\widehat X$, let
$\DD_\FF = \bigcup\{\CC_B : B \in \FF\}$.
The sets $\CC_B$ and $\DD_\FF$ are subgroups of $X$.
$\CC_B$ always has Haar measure 0, while the measure of 
$\DD_\FF$ depends on $\FF$.
We show that there is a filter $\FF$ such that
$\DD_\FF$ has measure $0$ but is not contained in any $\CC_B$.
This generalizes previous results for the special case
where $X$ is the circle group.
\end{abstract}

\section{Introduction} 
\label{sec-intro}
In this paper we study the pointwise convergence of sequences of characters
of compact abelian groups and its relation to Bohr topologies.
We begin with some abstract definitions.
All spaces considered here are assumed to be Hausdorff.

\begin{definition}
If $X,Y$ are topological spaces, then
$C(X,Y)$ is the set of continuous functions from $X$ to $Y$,
and $C_p(X,Y)$ denotes $C(X,Y)$
given the topology of \emph{pointwise convergence} \textup{(}i.e., regarding
$C_p(X,Y)$ as a subset of $Y^X$ with the usual product topology\textup{)}.
If $Y$ contains a distinguished point $1$, then
$\one$ denotes the constant function $x \mapsto 1$ in $C(X,Y)$.
\end{definition}

See Arkhangel'skii \cite{AR} for a discussion of such function spaces.

Suppose $X$ is a compact abelian group
and $Y = \TTT \subset \CCC$,
where $\TTT$ is the circle group.
As usual (see \cite{FOL,HR, RUD}),  $\duX$
denotes the dual group of $X$; that is, the group of characters,
or continuous homomorphisms into $\TTT$; then $\one$ is the
identity element of $\duX$.
If $G = \duX$ and we view $G$ as a discrete abelian group, then
$X \cong \widehat G$ by the Pontrjagin Duality Theorem.
However, if we consider $G \subseteq C_p(X,\TTT)$, then its inherited
topology is the \textit{Bohr topology} on $G$, and the closure
of $G$ in $\TTT^X$ is the \textit{Bohr compactification}, $\bohr G$, of $G$.
$G^\#$ denotes $G$ with its Bohr topology.
Since the compact group $\bohr G$ is dense in itself, and
$G^\#$ is dense in $\bohr G$, we have:

\begin{lemma}
\label{lemma-dense}
If $X$ is an infinite compact abelian group, then $\duX$
is dense in itself in the topology inherited from $C_p(X,\TTT)$.
\end{lemma}

However, $\duX$ has no pointwise convergent sequences.
To study pointwise convergence, we use the following notation:

\begin{definition}
\label{def-CB}
If $X,Y$ are topological spaces, $y \in Y$, and $B \subseteq C(X,Y)$
is infinite, then
$\CC_B(y)$ is the set of all $x \in X$ such that the
sequence $\langle \varphi(x) : \varphi \in B \rangle$ converges to $y$
\textup{(}that is, every neighborhood of $y$ contains $\varphi(x)$ for all
but finitely many $\varphi \in B$\textup{)}.
$\wCC_B = \bigcup_{y \in Y} \CC_B(y)$.
If $Y$ is a topological group with identity $1$,
then $\CC_B$ denotes $\CC_B(1)$.
\end{definition}

If $X$ and $Y$ are topological
groups and $B$ is a family of homomorphisms,
then $\CC_B$ and $\wCC_B$ are subgroups of $X$.
Clearly, $\CC_B \subseteq \wCC_B$.
The sequence $\langle \varphi : \varphi \in B \rangle$ converges
pointwise (i.e., in $C_p(X,Y)$) iff $\wCC_B = X$. 
So when $X$ is compact abelian and $B \subseteq \duX$,
$\; \wCC_B$ can never equal $X$,
but it can be non-trivial.  In \S\ref{sec-elem} we prove the following,
which gives some results involving the sizes of
$\CC_B$ and  $\wCC_B$:

\begin{theorem}
\label{thm-CB}
Let $X$ be  an infinite compact abelian group with $G = \duX$.
Then:
\begin{itemizz}
\item[1.] $\wCC_B$ is a Haar null set for each infinite $B \subseteq G$.
\item[2.] For any countable $Q \subseteq X$,
there is an infinite $B \subseteq \duX$ such that
$Q \subseteq \CC_B$, $\CC_B$  contains a perfect subset,
and $\CC_B$ is dense in $X$.
\item[3.] $\lambda(\overline B) \le 1/|\wCC_B|$ for all
infinite $B \subseteq G$.   Here, $\overline B$ is the closure
of $B$ in $\bohr G$, $\lambda$ is the Haar probability measure on $\bohr G$,
and $1/|\wCC_B| = 0$ when $|\wCC_B|$ is infinite.
\end{itemizz}
\end{theorem}

So, $\wCC_B$ is small in the sense of measure, but by (2),
even the smaller $\CC_B$ can be ``big'' in some senses.
However, (3) implies that whenever 
$\wCC_B$ is infinite, $B$ must be ``thin'' in $G$, in the sense that
$\overline B$ is a Haar null set in $\bohr G$.

When $X = \TTT$, the fact that
$\CC_B$ is null  is pointed out in \cite{BDMW, CTW}.

Note that both $\CC_B$ and $\wCC_B$
get bigger as $B$ gets smaller, so that the detailed arguments in
this paper will only involve countable $B$.  For example,
it is sufficient to prove (1) for countable $B$, and the
$B$ produced in the proof of (2) will be countable.

If $X = \TTT$ then $\widehat \TTT \subset C(\TTT,\TTT)$ is the
set of functions $z \mapsto z^n$ for $n\in\ZZZ$;
we identify $\widehat \TTT$ with $\ZZZ$.  As an illustration of (3),
let $B = \{kn : n \in \ZZZ\}$.  Then $\CC_B = \wCC_B = \{z\in \TTT : z^k = 1\}$,
and $\lambda(\overline B) = 1/k = 1/|\wCC_B| = 1/|\CC_B|$.

When $X = \TTT$, 
Barbieri,  Dikranjan,  Milan, and  Weber \cite{BDMW} showed
that assuming Martin's Axiom, there is a 
Haar null subgroup $D$ of $\TTT$ which is not contained in any $\CC_B$.
In \cite{HKcp} we showed that this holds in ZFC; in fact,
we gave an explicit definition of such a $D$ which is a Borel
set in $\TTT$.

There are two natural generalizations of these results about $C(\TTT,\TTT)$.
First, one may study the maps $(z \mapsto z^n) \in C(X,X)$ for
any compact group $X$; this was done in \cite{HKcp}.
In this paper, we consider the second generalization.  For
an arbitrary compact abelian group $X$,
we have $B \subseteq \duX \subset C(X,\TTT)$.
We shall produce (Theorem \ref{thm-main})
a Haar null subgroup $D$ of $X$ such that
$D$ is not contained in any countable union of the form
$\bigcup_\ell \wCC_{B_\ell}$.
As in \cite{HKcp}, it is convenient to define the null group $D$ 
from a filter:

\begin{definition}
Suppose that $X,Y$ are topological spaces, $y \in Y$, and
$\FF$ is a free filter on the set $C(X,Y)$.
Then $\DD_\FF(y) = \bigcup \{\CC_B(y) : B \in \FF\}$, and
$\wDD_\FF = \bigcup \{\wCC_B : B \in \FF\}$.
If $Y$ is a topological group with identity $1$,
then $\DD_\FF$ denotes $\DD_\FF(1)$.
\end{definition}

As usual, a filter $\FF$ is \textit{free} iff it contains the complements
of finite sets.  
As in \cite{HKcp}, our null group $D$ will be $\DD_\FF$, where
$\FF$ is a filter of sets of asymptotic density one:

\begin{definition}  
\label{def-dens}
For $E \subseteq \omega$, let $\underline{d}(E)$ and 
$\overline{d}(E)$ denote the lower and upper 
\emph{asymptotic density}:
\[
\underline{d}(E) =
\liminf_{n\to\infty} \frac{ |E \cap n|}{ n } \le
\limsup_{n\to\infty} \frac{ |E \cap n|}{ n } =
\overline{d}(E) \ \ .
\]
If equality holds, let 
$d(E)= \underline{d}(E) = \overline{d}(E)$
denote \emph{the asymptotic density} of $E$.
\end{definition}

\begin{definition}  
\label{def-filt}
Let $X$ be a compact abelian group and let 
$\vphi = \langle \varphi_n : n \in \omega \rangle$ be a
sequence of distinct elements of $\duX$.  Then
$\Fphi$ is
the filter $\FF$ generated by all sets of the form $\{\varphi_n : n \in E\}$
such that $d(E) = 1$.
\end{definition}  

\begin{proposition}
\label{prop-filt-null}
For $\Fphi$ defined as in \ref{def-filt},
$\wDD_{ \Fphi}$
is a Haar null subgroup of $X$.
\end{proposition}

Note that $\wDD_{ \Fphi}$ is clearly a subgroup.
We prove that it is null in \S\ref{sec-elem}.
The group $\wDD_{ \Fphi}$ could be trivial;
for example, if $X = \TTT$ and $\varphi_n(z) = z^n$, then
$\wDD_{ \Fphi} = \{1\}$.
In \cite{HKcp}, our null subgroup of $\TTT$ was of the form 
$\wDD_{ \Fphi}$, where $\varphi_n(z) = z^{n!}$.

The null group $\wDD_{\Fphi}$ contains $\DD_{\Fphi}$.
Nevertheless, Theorem \ref{thm-main} shows that for suitable $\vphi$,
even $\DD_{ \Fphi}$ is not contained in
any countable union of $\wCC_B$ sets.

\begin{theorem}
\label{thm-main}
For any  infinite compact abelian group $X$, there is a $D$ such that:
\begin{itemizz}
\item[1.]
$D$ is a Haar null subgroup of $X$;
\item[2.]
$D$ is dense in $X$;
\item[3.]
$D$ is not
a subset of any countable union of the form
$\bigcup_\ell \wCC_{B_\ell}$,
where each $B_\ell$ is an infinite subset of $\widehat X$;
\item[4.]
$D = \DD_{ \Fphi}$ for some sequence of distinct characters
$\vphi = \langle \varphi_n : n \in \omega \rangle$.
\end{itemizz}
\end{theorem}

The proof of Theorem \ref{thm-main} has two parts.  
In \S\ref{sec-nice}, we prove the theorem when $X$ is one of
four types of ``stock'' compact groups.
And  in \S\ref{sec-prelim}, we show
that it is sufficient to prove the theorem for those stock groups.
This argument applies the structure theory for abelian groups
to $\duX$, and is similar to the analysis used in constructing
$I_0$ sets (Hartman and Ryll-Nardzewski \cite{HAR}, Thm.~5;
see also \cite{KR}).

The stock groups are all second countable (that is, their
$\duX$ are countable).  The $|\duX|$ in Proposition
\ref{prop-filt-null} and Theorem \ref{thm-main} can be an arbitrary infinite
cardinal.  However, since $\wCC_B$ gets
bigger as $B$ gets smaller, it is sufficient to prove
Theorem \ref{thm-main} in the case that all the $B_\ell$ are countable.
For countable $B$,  $\CC_B$ and $\wCC_B$ are Borel 
(in fact, $F_{\sigma\delta}$) sets; likewise,
$\DD_{ \Fphi}$ and $\wDD_{ \Fphi}$ are $F_{\sigma\delta}$ sets
(see Proposition \ref{prop-borel}).

Our results are related to the notions of $\gggg$-closure 
and $\gggg$-density described by
Dikranjan, Milan, and Tonolo \cite{DMT}.
These notions may be expressed in terms of
an intersection involving our $\CC_B$:

\begin{definition}
Let $X$ be a compact abelian group, and
$J \le X$, with $\overline J$ its 
\textup{(}usual topological\textup{)} closure.
Then
$\gggg_X(J) = \overline J  \; \cap \;  \bigcap
\{\CC_B : B \in [\duX]^{\aleph_0} \ \&\ J \subseteq \CC_B \}$.
\end{definition}

\noindent
They call $\gggg_X(J)$ the \textit{$\gggg$-closure} of $J$ and
say that
$J$ is \textit{$\gggg$-dense} iff $\gggg_X(J) = X$.
Barbieri,  Dikranjan,  Milan, and  Weber \cite{BDMW2} ask
(see Question 5.7) whether for every infinite compact abelian group,
there is a $\gggg$-dense subgroup which is a Haar null set,
and they provide an affirmative under Martin's Axiom in some cases.
Our $D$ from Theorem \ref{thm-main} provides an affirmative answer
in all cases in ZFC.

\section{Elementary Facts}
\label{sec-elem}
Proposition \ref{prop-filt-null} is easily proved using
Ces\`aro limits:

\begin{definition}
\label{def-Cesaro-limit}
Given $r_n \in \CCC$ for $n \in \omega$ and $s \in \CCC$,
$r_n \leadsto s$ means that $\frac{1}{j} \sum_{n< j} r_n$ 
converges to $s$ as $j \to \infty$.
\end{definition}

\begin{lemma}  
\label{lemma-stat-Csearo}
Fix $r_n \in \CCC$ for $n \in \omega$ and $s \in \CCC$.
Assume that there is an $M \ge 0$ such that $|r_n| \le M$ for all $n$, and that 
${\lim_{n\in E}r_n = s}$ for some $E \subseteq \omega$ with $d(E) = 1$.
Then $r_n \leadsto s$.
\end{lemma}  

The following is proved exactly like Lemma 4.9
of  \cite{KR}, although the basic idea for the proof
goes back to Weyl \cite{WE}\S 7.

\begin{lemma}
\label{lemma-D-orthog}
Let $\mu$ be a probability measure on $X$.
Let $\varphi_n : X \to \CCC$, for $n \in \omega$, be measurable.
Assume that $M \ge 0$, 
$|\varphi_n(x)| \le M$ for all $n$ and $x$, and 
the $\varphi_n$ are orthogonal in $L^2(\mu)$.
Then $\mu(\{x \in X: \varphi_n(x) \leadsto 0\}) = 1$.
\end{lemma}

\begin{proofof}{Proposition \ref{prop-filt-null}}
Use Lemma \ref{lemma-stat-Csearo} and \ref{lemma-D-orthog}.
Here, the $\varphi_n$ map into $\TTT$, so
$\DD_{\Fphi}(0) = \emptyset$, so that
$\wDD_{ \Fphi}$ is disjoint from 
$\{x \in X: \varphi_n(x) \leadsto 0\}$.
\end{proofof}

The next lemma is immediate from the Pontrjagin Duality Theorem:

\begin{lemma}
\label{lem-contra}
For compact abelian $X$ and $Y$,
if $\widehat Y$ is isomorphic to a subgroup of $\widehat X$,
then there is a continuous homomorphism $\pi$
mapping $X$ onto $Y$.
\end{lemma}

Given compact abelian $X$, we can choose $Y$ so that $\widehat Y$ is a
countable subgroup of $\widehat X$.  Then $Y$ is second countable.
This sometimes lets us reduce a statement about arbitrary $X$
to a statement about second countable groups, as is illustrated
in the proof below of {Theorem \ref{thm-CB}(2).
It is also useful to recall:

\begin{lemma}
\label{lem-clopen}
If $\pi$ is a continuous homomorphism mapping the compact group
$X$ onto $Y$, then $\pi$ is both a closed map and an open map.
Also, $\lambda_X(\pi\iv(E)) = \lambda_Y(E)$ for all Haar-measurable
$E \subseteq Y$, where 
$\lambda_X,\lambda_Y$ are the Haar probability measures
on $X,Y$, respectively.
\end{lemma}

To prove  {Theorem \ref{thm-CB}(3), we need:

\begin{lemma}
\label{lem-bohr-null}
Every  infinite discrete abelian group $G$ is a Haar
null subset of $\bohr G$.
\end{lemma}

This lemma is immediate from
Varopoulos \cite{VAR}, who proves a more general result.
To prove the result directly for discrete abelian groups, note that
for countable ones, the result is trivial.
So for an arbitrary infinite discrete abelian $G$,
take a homomorphism $\pi$ from $G$ onto a countable $H$,
and then note that $\pi$ induces $\bohr \pi : \bohr G \onto \bohr H$,
with $G \subseteq (\bohr \pi)\iv (H)$.

The following lemma is also needed for  {Theorem \ref{thm-CB}(3):

\begin{lemma}
\label{lem-measure}
Let $X$ be a compact abelian group with $G = \duX$, and fix $u \in X$
and a subgroup $S$ of $G$.  Let
$K = \{x \in X : \forall \varphi \in S \, [\varphi(x) = \varphi(u)] \}$.
Then $\lambda(K) = 1 / |S|$, where $\lambda$ is Haar measure on $X$.
\end{lemma}
\begin{proof}
Let $\pi : X \onto \widehat S$ be the natural map.
Viewing $X$ as the characters of $G$, we have
\[
K = \{x \in X : \forall \varphi \in S \, [x(\varphi) = u(\varphi)] \} =
\{x \in X : x \res S = u \res S \}  = \pi\iv\{u\res S\}\ \ .
\]
Here, $u\res S$ is a point in $\widehat S$.
Since $\pi$ preserves Haar measure (see Lemma \ref{lem-clopen}),
if $S$ is infinite then $\lambda(K) = 0$, while if $S$ is
finite then $\lambda(K) = 1/|\widehat S|  = 1/|S|$.
\end{proof}

\begin{proofof}{Theorem \ref{thm-CB}}
Part (1) is clear from Proposition \ref{prop-filt-null},
since it is sufficient to prove it when $B$ is countable.

For (2), we shall produce a perfect subset of $X$ via a tree of
open sets indexed by finite 0-1 sequences.
List $Q$ as $\{q_j : j \in \omega\}$.
We now get distinct $\varphi_n \in \duX$ for $n \in \omega$ and 
$U_s \subseteq X$ for $s \in 2^{<\omega} =
\bigcup_{n \in \omega} \{0,1\}^n$
so that:
\begin{itemizn}{"2B}
\item $U_s$ is open and nonempty.
\item $\cl(U_{s\cat 0}) \cap \cl(U_{s\cat 1}) = \emptyset$ and
$\cl(U_{s\cat 0}),  \cl(U_{s\cat 1}) \subseteq U_s$.
\item $|1 - \varphi_n(x)| < 1/n$ whenever
$x \in \{q_j : j \le n\} \cup \bigcup \{U_s : s \in 2^n\}$.
\end{itemizn}
We do this by induction on $n$.  $\varphi_0$ can be arbitrary
and $U_{()}$ can be $X$.
If we are given
$U_s$ for $s \in 2^n$ and $\varphi_0, \ldots , \varphi_n$: 
First, choose distinct $p_{s\cat 0}, p_{s\cat 1} \in U_s$.
Then choose $\varphi_{n+1} \notin \{\varphi_0, \ldots , \varphi_n\}$
such that  $|1 - \varphi_{n+1}(x)| < 1/(n+1)$ whenever
$x \in \{q_j : j \le n+1\} \cup \{p_t : t \in 2^{n+1}\}$;
this is possible because $\one \in \duX \subset C_p(X,\TTT)$ and
is not isolated in $\duX$ (see Lemma \ref{lemma-dense}).
Then, we may choose $U_t$ for $t \in 2^{n+1}$ using 
the continuity of $\varphi_{n+1}$.

Let $K = \bigcap_{n\in\omega} \bigcup_{s\in 2^n} \cl(U_s)$,
and let $B = \{\varphi_n : n\in\omega\}$.
Then $K \cup Q \subseteq \CC_B$.  $K$ is not scattered, since
it maps continuously onto the Cantor set, so its perfect kernel
is non-empty.

We still need to get $\CC_B$ dense in $X$.  If $X$ is separable,
this is trivial, since 
we may assume that the countable $Q$ contains
a dense subset of $X$.
So for any $X$, choose a second countable $Y$
with $\widehat Y < \duX$, and
let $\pi : X \onto Y$ be as in Lemma \ref{lem-contra}.
For this separable $Y$, choose  an infinite $B_Y \subseteq \widehat Y$
such that $\pi(Q) \subseteq \CC_{B_Y}$,
$\CC_{B_Y}$ is dense in $Y$, and $\CC_{B_Y}$ contains
a perfect subset.
Then $B = \{\varphi \circ \pi : \varphi \in B_Y\}$
satisfies (2);
since $\pi$ is an open map (Lemma \ref{lem-clopen}),
$\CC_B = \pi\iv(\CC_{B_Y})$  is dense in $X$.

For (3), define $\Theta_0 : \wCC_B \to \TTT$ so that
$\Theta_0(x) $ is the limit of 
the sequence $\langle \varphi(x) : \varphi \in B \rangle$
(which exists by definition of $\wCC_B$).  Note that
$\Theta_0$ is a homomorphism from the group $\wCC_B$ into $\TTT$,
so, since $\TTT$ is divisible, it extends to a homomorphism
$\Theta : X \to \TTT$.  Then
$\wCC_B =
\{x \in X : \langle \varphi(x) : \varphi \in B \rangle \to \Theta(x)\}$.
Let $X_d$ denote the group $X$ with the discrete topology;
then we can identify $\bohr G$ with the compact group $\widehat{X_d}$.
So, $\Theta \in \bohr G$.  We can view $G^\#$  as a dense subgroup
of $\bohr G$, so that each $x \in X$ can be identified
with a continuous homomorphism on $\bohr G$.
With this identification, each $x \in \wCC_B$ satisfies
$\langle x(\varphi) : \varphi \in B \rangle \to x(\Theta)$, so that
$x(\Phi) = x(\Theta)$ for each $\Phi \in \overline B \backslash B$.
Thus, $\overline B \backslash B \subseteq 
\{\Phi \in \bohr G : \forall x \in \wCC_B \, [x(\Phi) = x(\Theta)] \}$,
so that $\lambda(\overline B \backslash B) \le 1 / |\wCC_B|$ by
applying Lemma \ref{lem-measure};
the $X,G,u,S$ in \ref{lem-measure} becomes
$\bohr G, X, \Theta, \wCC_B$ here.
Finally, $\lambda(\overline B) \le 1/|\wCC_B|$ because
$B \subseteq G$, which is a Haar null set in $\bohr G$ by
Lemma \ref{lem-bohr-null}.
\end{proofof}

Using the quotient argument in the last paragraph
in the proof of (2), getting $\CC_B$ to
contain a perfect set is trivial in ``most'' cases: if $G$
has any infinite subgroup with infinite index, then
$\CC_B$ will contain an infinite compact subgroup.

\section{Reduction to Stock }
\label{sec-prelim}
In this section, we show that it is sufficient to prove
Theorem \ref{thm-main} in the case that $X$ is the dual
of one of the groups listed in the following lemma:

\begin{lemma}
\label{lem-char-gps}
Every infinite abelian group contains a  subgroup isomorphic to
one of the following:
\begin{itemizn}{"2B}
\item $\ZZZ$.
\item $\sum_{n\in\omega} \ZZZ_{p_n}$, where the $p_n$ are primes
and $p_0 < p_1 < \cdots$.
\item $\sum_{n\in\omega} \ZZZ_{p}$, where $p$ is a fixed prime.
\item $\ZZZ_{p^\infty}$, for some prime $p$.
\end{itemizn}
\end{lemma}

This lemma is part of the structure theory for infinite abelian groups
(see Kaplansky \cite{KAP},
or Hewitt and Ross \cite{HR},
or \cite{KR} \S3).
The duals of these four groups are, respectively,
$\TTT$, $\prod_{n\in\omega} \ZZZ_{p_n}$, $(\ZZZ_p)^\omega$,
and the $p$-adic integers; for the last one, see \cite{HR} \S25.

Next, we use the $\pi : X \onto Y$ 
obtained from Lemma \ref{lem-contra} to translate a
$\vphi$ satisfying Theorem \ref{thm-main} for $Y$ 
to a $\vphi\circ\pi $ satisfying Theorem \ref{thm-main} for $X$.

\begin{lemma}
\label{lem-reduce}
Let $X$ and $Y$ be compact abelian groups,
with $\pi$ a continuous 
homomorphism  mapping $X$ onto $Y$.
Assume that 
$\vphi =  \langle \varphi_n : n \in \omega \rangle$ 
is a sequence of distinct characters of $Y$ such
that $\DD_{\Fphi}$ is not
a subset of any countable union of the form
$\bigcup_\ell \wCC_{A^\ell}$,
whenever each $A^\ell$ is an infinite subset of $\widehat Y$.
Let 
${\vphi\circ\pi} =
\langle \varphi_n\circ\pi : n \in \omega \rangle$.
Then, in $X$, 
$\DD_{\FF_{\vphi\circ\pi}  }$ is not
a subset of any countable union of the form
$\bigcup_\ell \wCC_{B^\ell}$,
whenever each $B^\ell$ is an infinite subset of $\widehat X$.
Also, if $\DD_{\Fphi}$ is dense in $Y$ then
$\DD_{\FF_{\vphi\circ\pi}  }$ is dense in $X$.
\end{lemma}

\begin{proof}
Let $K = \ker(\pi)$.
Since $\pi$ is an epimorphism, $X/K \cong Y$,
so characters of $Y$ correspond to characters of $X/K$.
Note also that each character in $\widehat X$
restricts to one in $\widehat K$.
Since $\wCC_B$ gets bigger as $B$ gets smaller,
we may shrink each $B^\ell$ to a countable set.
Shrinking  again to
$B^\ell = \{\psi^\ell_n : n \in \omega\}$,
we may assume that for each $\ell$, 
the $\psi^\ell_n \res K$, for $n \in \omega$, are either all the same 
or are all different.

\textit{Case 1}: 
The $\psi^\ell_n \res K$, for $n \in \omega$, are all the same. 
So each $\psi^\ell_n \cdot (\psi^\ell_0)\iv$
is identically $1$ on $K$,  and hence yields a character 
$\delta^\ell_n \in \widehat{X/K} \cong \widehat Y$
(with $\psi^\ell_n \cdot (\psi^\ell_0)\iv = \delta^\ell_n \circ \pi$).
Let $A^\ell = \{\delta^\ell_n : n \in \omega\} \subseteq \widehat Y$.
By our assumption on $\vphi$, we can fix a $y \in Y$
such that $y \in \DD_{\Fphi}$ and
$y \notin \wCC_{A^\ell}$ for all Case 1 $\ell$.
Note that if $x$ is \textit{any} element of $\pi\iv\{y\}$,
then $x \in \DD_{\FF_{\vphi\circ\pi}  }$.
Also, such an $x$ is not in $\wCC_{B^\ell}$ for all Case 1 $\ell$,
because 
the non-convergence of $\langle \delta^\ell_n(y) : n \in \omega \rangle$
implies the non-convergence of
$\langle \psi^\ell_n(x) \cdot (\psi^\ell_0(x))\iv  : n \in \omega \rangle$,
and hence the non-convergence of
$\langle \psi^\ell_n(x) : n \in \omega \rangle$.
We are thus done if we
produce $x\in \pi\iv\{y\}$ so that
$x \notin \wCC_{B^\ell}$ for all Case 2 $\ell$.
Fix $x^* \in \pi\iv\{y\}$.  Then our desired $x$ will be
an element of the coset $K x^* = \pi\iv\{y\}$.

\textit{Case 2}: The $\psi^\ell_n \res K$, for $n \in \omega$, are
all different.
For all Case 2 $\ell$, 
define $f^\ell_n : K \to \TTT$  by
$f^\ell_n(t) = \psi^\ell_n(t x^*) = \psi^\ell_n(t) \cdot \psi^\ell_n(x^*)$.
Note that each $f^\ell_n$
is the product of a character $\psi^\ell_n \res K$ of $K$
with a number $\psi^\ell_n(x^*)$,
so that $\{f^\ell_n : n \in \omega\}$ is an orthogonal 
family in $L^2(K)$.
It follows, by using Lemma \ref{lemma-D-orthog}, that 
$\widetilde{\CC}_{\{f^\ell_n : n \in \omega\}}$
is a Haar null set in $K$.
Choose $t$ such that
for each Case 2 $\ell$ the sequence
$ \langle f^\ell_n(t) : n \in \omega \rangle $ does not converge; then
$t x^* \notin \wCC_{B^\ell}$.

Finally, to prove  that 
$\DD_{\FF_{\vphi\circ\pi}  }$ is dense in $X$,
use the facts that $\pi$ is an open map by Lemma \ref{lem-clopen},
and that $\DD_{\FF_{\vphi\circ\pi}  } = \pi\iv(\DD_{\Fphi})$.
\end{proof}

\section{Nice Groups}
\label{sec-nice}
Definition \ref{def-nice} below isolates the key property of the
groups $\widehat G$, for the groups $G$ listed in
Lemma \ref{lem-char-gps}.

\begin{definition}
If $X$ is a compact abelian group, then
\[
\duX(X) = \{\varphi(x) : \varphi \in \duX \ \&\ x \in X\} \ \ .
\]
\end{definition}

\begin{proposition}
$\duX(X)$ is a subgroup of $\TTT$.
\end{proposition}
\begin{proof}
Let $G = \duX$.  If $G$ contains an element of infinite order,
then $\duX(X) = \TTT$.  Otherwise, $\duX(X)$ is the group generated by
all $e^{2\pi i / p^n}$ such that $p$ is prime and
$G$ contains an element of order $p^n$.
\end{proof}

If $G$ is of finite exponent (= bounded order), then 
$\duX(X)$ is finite; otherwise, $\cl(\duX(X)) = \TTT$.

\begin{definition}
\label{def-nice}
The compact abelian group $X$ is \emph{nice} iff $|\duX| = \aleph_0$ and
for all non-empty open $U \subseteq X$ and all $\varepsilon > 0$:
$\cl(\duX(X)) \subseteq N_\varepsilon(\varphi(U))$ for all
but finitely many $\varphi \in \duX$.
Here, $N_\varepsilon(S) =
\{z \in \TTT : \exists w \in S \,[ |z - w| < \varepsilon]\,\}$.
\end{definition}

\begin{lemma}
\label{lemma-nice-torsion}
If $G = \duX$ is an infinite torsion abelian group and
$\{\varphi\in G: \varphi^k = 1\}$ is finite for each $k$,
then $X$ is nice.
\end{lemma}
\begin{proof}
Note that $\cl(\duX(X)) = \TTT$, so we fix a
non-empty open $U \subseteq X$ and an $\varepsilon > 0$, and we must
verify that $N_\varepsilon(\varphi(U)) = \TTT$ for all
but finitely many $\varphi \in \duX$.
Observe that $\varphi(X)$ is finite for all $\varphi \in G$.
Translating $U$ and shrinking it, we may assume 
that $U = \{x \in X : \forall \psi \in F [ \psi(x) = 1 ]\}$,
where $F$ is a finite subgroup of $G$.
Let $R_n = \{z \in \TTT : z^n = 1\}$, and fix $m$
such that $N_\varepsilon(R_m) = \TTT$.
For all but finitely many $\varphi \in G$, the order of $[\varphi]$
in $G/F$ is at least $m$.  Fix any such $\varphi$; then for some $n \ge m$,
$\varphi^n \in F$ but $\varphi^k \notin F$ whenever $0 < k < n$.
Fix $y \in \Hom(G/F,\, \TTT)$ such that
$y([\varphi]) = e^{2 \pi i / n}$, this lifts to an
$x \in \widehat G = \Hom(G, \TTT)$ such that
$x(\varphi) =  e^{2 \pi i / n}$ and $x(\psi) = 1$ for all $\psi \in F$.
Identifying $\widehat G$ with $X$, we have $x\in U$ and
$\varphi(x) = e^{2 \pi i / n}$, so that $\varphi(U) \supseteq R_n$ because
$U$ is a group.  Then $n \ge m$ yields $N_\varepsilon(\varphi(U)) = \TTT$.
\end{proof}

\begin{lemma}
\label{lemma-nice-all}
$\widehat G$ is nice whenever $G$ is one of the groups listed in
Lemma \ref{lem-char-gps}.
\end{lemma}
\begin{proof}
Lemma \ref{lemma-nice-torsion} handles the duals of
$\sum_{n\in\omega} \ZZZ_{p_n}$ and $\ZZZ_{p^\infty}$.
For $\TTT = \widehat \ZZZ$,
note that for a given $U$,  $\varphi(U) = \TTT$ for all 
but finitely many $\varphi$.

For $G = \sum_{n\in\omega} \ZZZ_{p}$ and $\widehat G = (\ZZZ_p)^\omega$, 
follow the proof of Lemma \ref{lemma-nice-torsion}.
$U$ and $F$ are exactly the same.  Now, $\duX(X) = R_p$,
and $\varphi(U) = R_p$ for all $\varphi \notin F$.
\end{proof}

We now proceed to prove Theorem \ref {thm-main} for nice groups.

\begin{definition}
Let $X$ be a compact $2^{\mathit{nd}}$ countable abelian group with
metric $\rho$
and let $G = \duX$.  A \emph{nice partition} for $(X,\rho)$
is a sequence $\langle \Phi_j : j \in \omega \rangle$ such
that the $\Phi_j$ are finite disjoint nonempty sets whose union is $G$
and, if we set
\[
\rho_j(x,y) = \rho(x,y) + \sum\left\{ |\varphi(x) - \varphi(y)| :
\varphi \in \bigcup_{k \le j} \Phi_k \right\}\ \ ,
\]
then
for each $j$ and all $\varphi \in \bigcup_{k \ge {j+2}} \Phi_k $,
all $x \in X$, and all $z \in \cl(\duX(X))$,
there is a $y \in X$ with $\rho_j(x,y) < 2^{-j}$
and $|\varphi(y) - z| < 2^{-j}$.
\end{definition}

\begin{lemma}
If $X$ is a nice compact abelian group with metric $\rho$,
then there is a nice partition for $(X,\rho)$.
\end{lemma}
\begin{proof}
List $G$ as $\{\varphi_j : j \in \omega\}$.  Now, define the $\Phi_j$
by induction.  Let $\Phi_0 = \{\varphi_0\}$.  Given
$\Phi_k$ for $k \le j$, we have the metric $\rho_j$ on $X$, so
for some finite $m$, we 
may cover $X$ by open sets $U_0, \ldots, U_m$ of $\rho_j$--diameter
less than $\varepsilon := 2^{-j}$.  Now choose $\Phi_{j+1}$ 
so that $\cl(\duX(X)) \subseteq N_\varepsilon(\varphi(U_\ell))$ for all $\ell$
and for all $\varphi \in G \setminus \bigcup_{k \le {j+1}} \Phi_k $.
Also make sure that $\varphi_j \in \bigcup_{k \le {j+1}} \Phi_k$  so that
$G$ will be the union of all the $\Phi_j$.
\end{proof}

\begin{definition}
Suppose that $\bm{\Phi} = \langle \Phi_j : j \in \omega \rangle$
is a nice partition for $(X,\rho)$.  A sequence 
$\langle \varphi_n : n \in \omega\rangle$ from $\duX$ is
\emph{thin} \textup{(}with respect to $\bm{\Phi}$\textup{)} iff each
$\varphi_n \in \Phi_{j_n}$, where each $j_{n+1} \ge j_n + 2$.
\end{definition}

\begin{lemma}
\label{lemma-split}
Assume that $\langle \varphi_n : n \in \omega\rangle$ is a thin sequence,
$\omega$ is partitioned into 
two infinite sets, $A,B$, and $a,b \in \cl(\duX(X))$.
Then for some $x\in X$,
\[
\varphi_n(x)
\raisebox{-0.7ex}{\mbox{$\, \overrightarrow{ \mbox{\tiny $n\in A$}\,}$ }} a
\qquad \mathrm{AND} \qquad
\varphi_n(x)
\raisebox{-0.7ex}{\mbox{$\, \overrightarrow{ \mbox{\tiny $n\in B$}\,}$ }} b
\ \ .
\]
\end{lemma}
\begin{proof}
Choose $x_n \in X$ for $n \in \omega$ as follows:  $x_0,x_1$ are arbitrary.
Given $x_{n-1}$ with $n\ge 2$, use 
$\varphi_n \in \bigcup_{k \ge {j_{n-1}+2}} \Phi_k $,
plus $j_{n-1} \ge 2(n -1) \ge n$, to get $x_{n}$ to satisfy:
\begin{itemizn}{"2B}
\item $\rho_{j_{n-1}}(x_{n-1}, x_n) < 2^{-n}$.
\item $n \in A \ \Rightarrow\  |\varphi_n(x_{n}) - a| < 2^{-n}$.
\item $n \in B \ \Rightarrow\  |\varphi_n(x_{n}) - b| < 2^{-n}$.
\end{itemizn}
\noindent
Then each $\rho(x_{n-1}, x_n) < 2^{-n}$, so
$\langle x_n : n \in \omega \rangle$ converges to some $x$.
Now, fix $n \ge 1$, and we estimate $|\varphi_n(x_n) - \varphi_n(x)|$:
For all $m > n$, 
$|\varphi_n(x_{m-1}) - \varphi_n(x_m)| \le
\rho_{j_{m-1}}(x_{m-1}, x_m) < 2^{-m}$.
Thus, $|\varphi_n(x_n) - \varphi_n(x)| \le
\sum_{m = n+1}^\infty2^{-m} = 2^{-n}$.

Now, if $n \in A$, then 
\[
|\varphi_n(x) - a| \le |\varphi_n(x_n) - \varphi_n(x)| +
|\varphi_n(x_{n}) - a|  \le 2^{-n} + 2^{-n} \to 0 \ \ .
\]
The argument is the same for $n \in B$.
\end{proof}

\begin{lemma}
\label{lemma-nice-are-nice}
Let $X$ be a compact $2^{\mathit{nd}}$ countable abelian group with metric
$\rho$
and let $G = \duX$.  Suppose that
$\bm{\Phi} = \langle \Phi_j : j \in \omega \rangle$
is a nice partition for $(X,\rho)$ and 
$\vphi = \langle \varphi_n : n \in \omega\rangle$,
where each $\varphi_n \in \Phi_{3n}$.
Let $B_\ell$, for $\ell \in \omega$, be any infinite subsets of $G$.
Then $\DD_{\Fphi} \not\subseteq \bigcup_\ell \widetilde{\CC}_{B_\ell}$.
\end{lemma}
\begin{proof}
By a standard diagonal argument, get $\varphi'_n$ for $n \in \omega$
and $E,F \subseteq \omega$ such that:
\begin{itemizz}
\item[1.] $\langle \varphi'_n : n \in \omega\rangle$ 
is thin with respect to $\bm{\Phi}$.
\item[2.] $d (\{n : \varphi'_n = \varphi_n \}) = 1$.
\item[3.] $E \cup F = \omega$ and $E \cap F = \emptyset$.
\item[4.] $d (E) = 1$.
\item[5.] For each $\ell$, both
$\{n \in E : \varphi'_n \in B_\ell\}$ and
$\{n \in F : \varphi'_n \in B_\ell\}$ are infinite.
\end{itemizz}
Fix $z \in \duX(X) \backslash \{1\}$. 
By (1)(3), we may apply Lemma \ref{lemma-split} and fix $x \in X$ such that
\[
\varphi'_n(x)
\raisebox{-0.7ex}{\mbox{$\, \overrightarrow{ \mbox{\tiny $n\in E$}\,}$ }} 1
\qquad \mathrm{AND} \qquad
\varphi'_n(x)
\raisebox{-0.7ex}{\mbox{$\, \overrightarrow{ \mbox{\tiny $n\in F$}\,}$ }} z
\ \ .
\]
By (2)(4), $x \in \DD_{\Fphi} $.
By (5), $x \notin  \widetilde{\CC}_{B_\ell}$ for each $\ell$.
\end{proof}

\begin{lemma}
\label{lemma-nice-dense}
Suppose that $\bm{\Phi} = \langle \Phi_j : j \in \omega \rangle$
is a nice partition for $(X,\rho)$ and
$\vphi = \langle \varphi_n : n \in \omega\rangle$ from $\duX$ is
thin with respect to $\bm{\Phi}$.  Let
$B = \{\varphi_n : n \in \omega\}$.  Then $\CC_B$ is dense in $X$,
so that $\DD_{\Fphi}$ is dense in $X$.
\end{lemma}
\begin{proof}
This is similar to the proof of Lemma \ref{lemma-split}.
Fix a non-empty open $U \subseteq X$.  We must
produce an $x \in U$ such that $\varphi_n(x) \to 1$.
We may assume that $q \in X$ and $r\in\omega$
and $U = \{x \in X : \rho(x,q) < 2^{-r+1}\}$.
Choose $x_n \in X$ for $n \in \omega$ as follows: 
$x_0 = x_1 = \cdots = x_r = q$.
Given $x_{n-1}$ with $n\ge r+1$,
get $x_{n}$ to satisfy:
\begin{itemizn}{"2B}
\item $\rho_{j_{n-1}}(x_{n-1}, x_n) < 2^{-n}$.
\item $|\varphi_n(x_{n}) - 1| < 2^{-n}$.
\end{itemizn}
\noindent
Then $\langle x_n : n \in \omega \rangle$ converges to some $x$
with $\rho(x,q) \le 2^{-r}$, so $x \in U$.
As in the proof of Lemma \ref{lemma-split}, $|\varphi_n(x) - 1| \to 0$.
\end{proof}

\begin{proofof}{Theorem \ref {thm-main}}
By Lemmas \ref{lemma-nice-are-nice} and \ref{lemma-nice-dense},
the theorem holds for all
nice groups, which by Lemma \ref{lemma-nice-all}, 
includes the duals of all the
groups listed in Lemma \ref{lem-char-gps}.
Then, by Lemma \ref{lem-reduce},  the theorem holds for all $X$.
\end{proofof}

Note that not every $X$ with a countable dual is nice;
see Example \ref{ex-not-nice}.

\section{Remarks and Examples}
\label{sec-rem}
The proof in \S\ref{sec-elem}
that $\wDD_{\Fphi}$ is null
makes essential use of asymptotic density, via
Lemma \ref{lemma-stat-Csearo}; one cannot replace
$\Fphi$ by an arbitrary filter $\FF$, 
since  $\wDD_\FF$, or even the smaller
$\DD_\FF$, might be all of $X$. 
By Proposition 1.2 of \cite{HKcp}  and
Lemma \ref{lemma-dense},

\begin{proposition}
\label{prop-filter}
If $X$ is any infinite compact abelian group, then there is
a free filter $\FF$ on $\duX$ such that 
$\FF$ contains a countable set and $\DD_\FF = X$.
\end{proposition}

It is not clear whether the nice groups are of interest in their own
right, or just an artifact in the proof of Theorem \ref{thm-main}.
Not every dual of a countable discrete abelian group is nice:
 
\begin{example}
\label{ex-not-nice}
$\ZZZ_4 \times (\ZZZ_2)^\omega$ is not nice.
\end{example}
\begin{proof}
Elements of $X = \ZZZ_4 \times (\ZZZ_2)^\omega$ are of the form
$\langle x, \vec y\rangle$, where $x \in \ZZZ_4 = \{1, i, -1, -i\}$,
$\ZZZ_2 = \{\pm 1 \}$, and
$\vec y \in (\ZZZ_2)^\omega$.
$\duX(X) = \{1, i, -1, -i\}$.
Let $\varphi_n(x, \vec y) = x \cdot y_n$.
Let $U = \{\langle x, \vec y\rangle : x = i\}$.
Then the $\varphi_n$ are distinct characters, and 
$\varphi_n(U) = \{\pm i\}$, so Definition \ref{def-nice} fails
whenever $\varepsilon < \sqrt{2}$.
\end{proof}

All the ``$\CC$'' and ``$\DD$'' sets discussed in this paper are Borel:

\begin{proposition}
\label{prop-borel}
Let $X$ be any compact abelian group.
If $B \subseteq \duX$ is
countably infinite, 
then $\CC_B$ and $\wCC_B$ are $F_{\sigma\delta}$ sets.
If $\vphi = \langle \varphi_n : n \in \omega \rangle$ is
a sequence of distinct elements of $\duX$, then
$\DD_{ \Fphi}$ and $\wDD_{ \Fphi}$ are $F_{\sigma\delta}$ sets.
\end{proposition}
\begin{proof}
Let $B = \{\varphi_n : n \in \omega\}$. 
Then $x \in \wCC_B $ iff
\[
\forall r \in \omega 
\, \exists s < r 
\, \exists k \in \omega 
\, \forall m > k
\left[ |\varphi_m(x) - e^{2\pi i s/r} | \le \frac{\pi}{r} \right] \ \ ,
\]
since $\TTT \subseteq \bigcup_{s < r} N_{\pi/r}(e^{2\pi i s/r} )$.
This displays
$\wCC_B$ as a countable intersection of $F_\sigma$ sets.
The argument for $\CC_B$ is similar; just replace $s$ by $0$.
Likewise, $x \in \wDD_{ \Fphi}$ iff
\[
\forall r \in \omega 
\, \exists s < r 
\, \exists k \in \omega 
\, \forall n > k
\left[
\frac{1}{n}
\left| \{m < n : |\varphi_m(x) - e^{2\pi i s/r} | \le \frac{\pi}{r} \} \right|
\ge 1 - \frac{1}{r} \right] \ \ .
\]
Again, replace  $s$ by $0$ to see that $\DD_{ \Fphi}$
is an $F_{\sigma\delta}$ set.
\end{proof}

It is natural to ask whether the countable $\{B_\ell : \ell\in \omega\}$
from Theorem \ref{thm-main} could be replaced by 
a family of $\aleph_1$ sets.
Under CH, this is clearly false, since then $|X|$ may be $\aleph_1$,
in which case Theorem \ref{thm-CB} implies that
a union of the form
$\bigcup_{\alpha < \omega_1} \CC_{B_\alpha}$ can be all of $X$.
Assuming Martin's Axiom (MA), the proof of Theorem \ref{thm-main}
implies that our $\DD_{ \Fphi}$ is not
a subset of any union of the form
$\bigcup_{\alpha < \kappa} \wCC_{B_\alpha}$ where
$\kappa < 2^{\aleph_0}$.  To see this, note that the countability of
the family $\{B_\ell : \ell\in \omega\}$ was only used in two places.
First, in handling Case 2 of Lemma \ref{lem-reduce}, we used
the fact that a compact group $K$ is not covered by $\aleph_0$ null sets,
and MA lets us replace the ``$\aleph_0$'' by ``$< 2^{\aleph_0}$''.
Second, the diagonal argument in the  proof Lemma \ref{lemma-nice-are-nice}
will work with families of size less than $2^{\aleph_0}$ under MA.

It is also consistent with ZFC to have $2^{\aleph_0}$ arbitrarily large
but $\TTT = \bigcup_{\alpha < \omega_1} \CC_{B_\alpha}$.
This proof resembles the standard construction of an ultrafilter
of character $\aleph_1$ (see \cite{KUNT}, Exercise VIII.A10).
Start with $2^{\aleph_0}$ large in the ground model $V$ and iterate
forcing $\aleph_1$ times with finite supports, forming $V_\alpha$ for
$\alpha \le \omega_1$.  When $\alpha < \omega_1$, 
let $\FF_\alpha \in V_\alpha$ be a filter on $\ZZZ = \widehat \TTT$
obtained from Proposition \ref{prop-filter}, and get
$B_\alpha \in V_{\alpha+1}$ so that
$B_\alpha \subset^* A$ for all $A \in \FF_\alpha$.
One can even make the $B_\alpha$ generate a P-point ultrafilter,
so that in the final model $V_{\omega_1}$, 
the $\FF$ of Proposition \ref{prop-filter} could be a P-point
of character $\aleph_1$.  To do this, make sure that each $B_\alpha$
is chosen so that $0$ is a limit point of $B_\alpha$
in the Bohr topology of $\ZZZ$.
Note that the $\FF$ of Proposition \ref{prop-filter} can never
be a selective ultrafilter, since it would then contain thin
sets and run afoul of Lemma \ref{lemma-split}.

\end{document}